\documentclass[referee]{gji}
\usepackage{timet,color}
\usepackage[urlcolor=blue,citecolor=black,linkcolor=black]{hyperref}

\usepackage[fleqn]{amsmath}
\usepackage{gji_extra}
\usepackage{graphicx}
\usepackage{placeins}

\title[Integral equation method with domain decomposition]
  {3D induction log modelling with integral equation method and domain decomposition preconditioning}
\author[D.H. Saputera et. al.]
  {\Large D.H. Saputera$^1$, M. Jakobsen$^1$, K.W.A. van Dongen$^2$, N. Jahani$^3$, K.S. Eikrem$^3$, S. Alyaev$^3$\\
  $^1$Department of Earth Science, University of Bergen, Bergen, Norway\\
  $^2$Department of Imaging Physics, Delft University of Technology, Delft, The Netherlands\\
  $^3$NORCE Norwegian Research Centre, Bergen, Norway\\
  }
\date{}
\pagerange{\pageref{firstpage}--\pageref{lastpage}}

\begin{document}

\label{firstpage}

\maketitle

\begin{summary}
The deployment of electromagnetic (EM) induction tools while drilling is one of the standard routines for assisting the geosteering decision-making process. The conductivity distribution obtained through the inversion of the EM induction log can provide important information about the geological structure around the borehole. To image the 3D geological structure in the subsurface, 3D inversion of the EM induction log is required. Because the inversion process is mainly dependent on forward modelling, the use of fast and accurate forward modelling is essential. In this paper, we present an improved version of the integral equation (IE) based modelling technique for general anisotropic media with domain decomposition preconditioning. The discretised IE after domain decomposition equals a fixed-point equation that is solved iteratively with either the block Gauss-Seidel or Jacobi preconditioning. Within each iteration, the inverse of the block matrix is computed using a Krylov subspace method instead of a direct solver. An additional reduction in computational time is obtained by using an adaptive relative residual stopping criterion in the iterative solver. Numerical experiments show a maximum reduction in computational time of 35 per cent compared to solving the full-domain IE with a conventional GMRES solver. Additionally, the reduction of memory requirement for covering a large area of the induction tool sensitivity enables acceleration with limited GPU memory. Hence, we conclude that the domain decomposition method is improving the efficiency of the IE method by reducing the computation time and memory requirement.
\end{summary}


\section{Introduction} \label{section_intro}

State-of-the-art tools for electromagnetic (EM) induction logging-while-drilling (LWD) enable real-time mapping of formation boundaries tens of meters away from the borehole \citep{sinha_past_2022}. These tools typically consist of multiple antenna configurations that have different sensitivities to the electrical resistivity distribution in the medium around the borehole. The distribution of the electrical properties is quantified through an inversion process and provides structural information and characteristics of the surrounding medium. The studies in real-time geosteering inversion usually employ 1D or 2D approximations \citep{bakr_fast_2017,noh_real-time_2022,pardo_fast_2015,puzyrev_deep_2019}. However, for imaging complex geological structures, it is important to capture the 3D variability of the resistivity change around the borehole through 3D inversion methods \citep{puzyrev_recent_2019,sinha_past_2022}. The work of \citet{wilson_real-time_2019} shows that it is possible to perform 3D inversion in real-time, however, it is challenging due to the large computational cost required for the 3D forward modelling, especially when quantification of the uncertainties in the inversion is required. Therefore, the study of a fast 3D forward solver that accurately models induction logs remains essential for the development and testing of new imaging methods.

The integral equation (IE) method is one of the most widely applied numerical methods for the 3D modelling of EM data \citep{avdeev_three-dimensional_2005,wang_geophysical_2021} alongside the finite difference \citep{newman2002three,hou2006finite} and finite element methods \citep{puzyrev2013parallel,ren2014hybrid}. One of the main advantages of using the IE method is that it has the accuracy of a semi-analytical solution \citep{wang_geophysical_2021}. Without introducing many specific approximations, the EM fields around the borehole are obtained by solving the linear system arising from the discretization of the integral equations. As the linear system is dense, the computational memory and time required can be large compared to other numerical methods \citep{yoon_hybrid_2016,zaslavsky_hybrid_2011}. To overcome this challenge, the linear system can be efficiently solved using an iterative solver based on the Krylov subspace method in combination with the utilization of the FFTs to accelerate the convolution integral operations in the linear system \citep{fang_efficient_2006}. A faster convergence rate can be achieved by implementing the contraction IE formulation \citep{hursan_contraction_2002} which works especially well in the presence of a high contrast or a high degree of anisotropy. Additionally, the application of GPUs further decreases computation times because GPUs enable the acceleration of mathematical operations that can be straightforwardly parallelized \citep{dyatlov_efficient_2015,saputera2022gpu}.

In the work of \citet{zhdanov_integral_2006}, the formulation of the IE method is extended by decomposing the region of interest into several sub-domains. The field in the entire domain is obtained by sequentially solving the linear system in each sub-domain and updating the interaction between the sub-domains iteratively until convergence. With this formulation, it becomes feasible to conduct large-scale modelling of surface EM data in heterogeneous media as the computational operation can be reduced to one sub-domain at a time. It is possible to obtain an additional reduction in computational costs by only considering sub-domains that contain an anomaly with respect to the background medium. This leads to a smaller number of discretization blocks required for the 3D modelling while still enabling FFT implementation \citep{endo_large-scale_2009} and an improved iterative solver convergence rate \citep{van_dongen_reduced_2007}. Typically, a horizontally layered model is chosen as the background medium as the theory of Green's functions for layered one-dimensional (1D) models is very well developed \citep{zhdanov_integral_2006}. Hence, the IE method can be very efficient when the resistivity model only deviates in some areas from the 1D model. However, in our application, the subsurface structure can vary in all directions. The sub-domains containing an anomaly can be everywhere around the EM tools and it may not be possible to achieve a reduction in the number of discretizations by the domain decomposition. Additionally, the sub-domains from the decomposition can be adjacent to each other such that the interactions between neighbouring sub-domains are not negligible. 

The domain decomposition method can lead to an efficient way of solving the linear system of the IE method \citep{jakobsen_distorted_2018,wang2017accurate}. In the work of \citet{jakobsen_distorted_2018}, the domain decomposition method is used to efficiently compute the T-matrix for the inversion of controlled source EM data. It is also shown that the domain decomposition method opens up the possibility to compute the T-matrix in parallel. 

In this paper, we demonstrate that the formulation of IE with domain decomposition (IE-DD) can be interpreted as a preconditioned linear system, offering a computational advantage. We illustrate that the IE-DD method can be represented as a fixed-point equation, which is iteratively solved using block Gauss-Seidel or Jacobi preconditioners \citep{saad_iterative_2003}. In particular, we will use a Krylov subspace method to invert the block matrices that are present in the formulation. Instead of expressing the decomposition formulation in terms of the contrast source in each sub-domain as described in \citet{zhdanov_integral_2006} and \citet{endo_large-scale_2009}, we present the domain decomposition formulation in terms of the electric field in each sub-domain and a different perspective on the derivation of the IE-DD formulation. Additionally, we propose the use of an inexact iterative solver when solving the IE linear system for each sub-domain where the target tolerance is adapted based on the full domain residual.

The outline of this paper is described as follows. In section \ref{section_theory} called theory, we give an overview of the theory and implementation of the conventional IE method and the IE-DD. In section \ref{section_numerical_results} called numerical results and discussion, we present three numerical examples to show the performance of the IE-DD method and discuss the computational aspect of our implementation. First, we show an example with isolated sub-domains to verify if the domain decomposition formulation will produce the same numerical results as the conventional full-domain formulation. In the second example, we show a numerical experiment with three different IE-DD schemes and compare the performance of these schemes with each other and the full-domain IE as a reference. In the last example, we simulate a logging scenario across a large complex 3D model. Furthermore, we showcase the ability of the domain decomposition method to reduce the memory requirement for dealing with a large number of grid blocks in the last example. This feature lets us cover more portion of the subsurface receivers while keeping a fine grid size, which may not be straightforward to implement in our currently available computer without the domain decomposition method. In section \ref{section_conclusion}, we provide a compact evaluation of the IE-DD implementation in this study and also some possible improvements for future research. This paper contains appendices with more in-depth details of the IE-DD derivation and implementation. We also include the comparison of our conventional IE code and existing code as a benchmark of our work in the appendix.

\section{Theory} \label{section_theory}

\subsection{The Integral Equation Method for 3D Induction Logs Modelling} \label{subsection_IE3D}

Maxwell's equations for heterogeneous media \citep{wannamaker2002three} are the basic theory for modelling the induction tools\textquotesingle{} response within the frequency domain:
\begin{gather}
    \nabla \times \mitbf{E}\left( \mitbf{r} \right) = i\omega\mu\mitbf{H}\left( \mitbf{r} \right) + \mitbf{J}^{H}\left( \mitbf{r} \right), \label{eq:E_maxwell}\\
    \nabla \times \mitbf{H}\left( \mitbf{r} \right) = \widehat{\mitbf{\sigma}}\left( \mitbf{r} \right)\mitbf{E}\left( \mitbf{r} \right), \label{eq:H_maxwell} 
\end{gather}
where \(\mitbf{E}\left( \mitbf{r} \right)\) and
\(\mitbf{H}\left( \mitbf{r} \right)\) are the total electric and magnetic fields, respectively, at location \(\mitbf{r}\),
\(\mitbf{J}^{H}(\mitbf{r})\) denotes the magnetic source term, $\omega$ is the angular frequency, \(\mu\) is the magnetic permeability, \(\widehat{\mitbf{\sigma}}\left( \mitbf{r} \right) = \mitbf{\sigma}\left( \mitbf{r} \right) - i\omega\mitbf{\varepsilon}\left( \mitbf{r} \right)\) is the complex electric conductivity, \(\mitbf{\varepsilon}\) is the dielectric permittivity, and $i$ = \(\sqrt{- 1}\). We assume that the magnetic permeability is constant and it is set equal to the magnetic permeability of the vacuum \(\mu_{0}\). Additionally, the imaginary part of the complex conductivity can be ignored in the diffusion regime, which is a typical assumption for the operating conditions of induction tools.

The total electric and magnetic fields can be formulated using the following integral equations \citep{fang_efficient_2006} 
\begin{gather}
    \mitbf{E}\left( \mitbf{r} \right) = \mitbf{E}^{(0)}\left( \mitbf{r} \right) + \int_{\Omega}^{}{\mitbf{G}^{E}\ \left( \mitbf{r},\mitbf{r'} \right)}\Delta\mitbf{\sigma}\left( \mitbf{r'} \right)\mitbf{E}\left( \mitbf{r'} \right)\,\, dV(\mitbf{r'}),\label{eq:E_IE}\\
    \mitbf{H}\left( \mitbf{r} \right) = \mitbf{H}^{(0)}\left( \mitbf{r} \right) + \int_{\Omega}^{}{\mitbf{G}^{H}\ \left( \mitbf{r},\mitbf{r'} \right)}\Delta\mitbf{\sigma}\left( \mitbf{r'} \right)\mitbf{E}\left( \mitbf{r'} \right)\,\, dV(\mitbf{r'}),\label{eq:H_IE}
\end{gather}
where the \(\Omega\) indicates the domain of integration where anomalies in the conductivity relative to the background conductivity are present. The integral terms in equations (\ref{eq:E_IE}) and (\ref{eq:H_IE}) represent the scattered electric and magnetic fields, respectively, due to the presence of these anomalies. The (0) superscripts indicate the fields defined for a homogeneous isotropic background medium with conductivity \(\sigma_{0}\) which are referred to as the background fields. We choose a homogeneous isotropic background medium for simplicity and efficiency when calculating Green's tensor \citep{fang_efficient_2006}, and we assume that the tool is not always surrounded by a horizontally layered medium. The tensor \(\Delta\mitbf{\sigma}\left( \mitbf{r} \right) = \mitbf{\sigma}\left( \mitbf{r} \right) - \sigma_{0}\mitbf{I}\), denotes the conductivity contrast between the actual anisotropic medium and the background medium, and with \(\mitbf{I}\) the identity tensor.
The electric Green's tensor \(\mitbf{G}^{E}\left(\mitbf{r,r'}\right)\) and its relation to the magnetic Green's tensor \(\mitbf{G}^{H}\mitbf{(r,r')}\) for a homogenous isotropic medium are \citep{fang_efficient_2006}
\begin{gather}
    \mitbf{G}^{E}\left( \mitbf{r},\mitbf{r'} \right) = \left( i\omega\mu_{0}\ \mitbf{I} + \frac{\nabla\nabla}{\sigma_{0}} \right)g\left( \mitbf{r},\mitbf{r'} \right),\label{eq:Green_E}\\
    \mitbf{G}^{H}\left( \mitbf{r},\mitbf{r'} \right) = \left( i\omega\mu_{0} \right)^{- 1}\nabla \times \mitbf{G}^{E},\label{eq:Green_H}\\
    g\left( \mitbf{r},\mitbf{r'} \right) = \frac{e^{ik_{0}\left| \mitbf{r} - \mitbf{r'} \right|}}{4\pi\left| \mitbf{r} - \mitbf{r'} \right|},\label{eq:Green_Scalar}
\end{gather}
where \(g\left( \mitbf{r},\mitbf{r'} \right)\) is the scalar Green's function and \(k_{0} = \sqrt{i\omega\mu_{0}\sigma_{0}}\). To calculate the total magnetic fields, the total electric fields need to be obtained first by solving equation (\ref{eq:E_IE}). Afterward, the calculation of the total magnetic fields is a straightforward addition of the background magnetic fields and the integral term as shown in equation (\ref{eq:H_IE}). Therefore, the main computational challenge of the IE method is to solve the integral equation (\ref{eq:E_IE}), which is classified as the Fredholm integral equation of the second kind \citep{fang_efficient_2006}.

\subsection{Numerical Implementation of the Integral Equation Method} \label{subsection_IE_numerics}

A numerical solution of the volume integral in equation (\ref{eq:E_IE}) can be obtained using the method of moments \citep{gibson_method_2021}. The subsurface model around the induction tool is discretized into a set of grid blocks with centroids \(\mitbf{r}^{j}\) and volume of \(\Delta v\), where $j$ indicates the $j$-th grid block. The discretization of equation (\ref{eq:E_IE}) leads to a linear system of equations that can be expressed in operator form as
\begin{equation} \label{eq:IE_linear_system}
    \left( \mitbf{I} - \mitbf{\mathcal{G}}\Delta\mitbf{\sigma} \right)\mitbf{E} = \mitbf{E}^{(0)},
\end{equation}
where \(\mitbf{\mathcal{G}}\) is the operator that represents the discrete convolution integral of the electric Green's tensor
\(\mitbf{G}^{E}\left( \mitbf{r},\mitbf{r'} \right)\) with the contrast source \(\Delta\mitbf{\sigma E}\) in equation (\ref{eq:E_IE}). The Green's function in equation (\ref{eq:Green_E}) can be discretized by separating the non-singular part of the Green's function and dealing with the singularity by integrating the Green's function of a single grid block over a spherical domain with an equivalent volume \citep{gao_analytical_2005,jakobsen_distorted_2018}. The linear system in equation (\ref{eq:IE_linear_system}) can be efficiently solved using a Krylov subspace method because it does not require the matrix of the linear system to be formed explicitly. The desired accuracy of the iterative method is quantified by the relative residual $e$ which is calculated as
\begin{equation} \label{eq:relres_definition}
    e = \frac{\left\| \mitbf{E}^{(0)} - \left( \mitbf{I} - \mitbf{\mathcal{G}}\Delta\mitbf{\sigma} \right)\mitbf{E} \right\|}{\left\| \mitbf{E}^{(0)} \right\|},
\end{equation}
where \(\left\| \cdot \right\|\) is the L\textsubscript{2}-norm. In this study, we use the generalized minimum residual or GMRES \citep{saad_gmres_1986} as the linear system solver.

Green's tensor operator exhibits a convolution structure in each of the tensor components. This property enables the use of FFT to convolve a Green's tensor component \(G_{pq}^{E}\) and a component of the contrast source \(\left( \Delta\mitbf{\sigma E} \right)_{q}\) efficiently \citep{fang_efficient_2006}. The $p$ and $q$ indices indicate the component of Green's tensor and the contrast source vector with $p$ and $q$ = x,y,z. At each step of the iterative solver, the convolution integral can be efficiently calculated by
\begin{equation} \label{eq:FFT_convolution}
    \mathcal{G}_{pq}\left( \Delta\mitbf{\sigma E} \right)_{q} = \mathcal{F}^{- 1}\left( \mathcal{F}\left\lbrack G_{pq}^{E} \right\rbrack\mathcal{\odot F}\left\lbrack \left( \Delta\mitbf{\sigma E} \right)_{q} \right\rbrack \right),
\end{equation}
where \(\mathcal{F}\) is the FFT operator and \(\odot\) denotes elementwise multiplication. This operation reduces the convolution computation complexity from O(\emph{N}\textsuperscript{2}) to O(\emph{N}log\textsubscript{2}\emph{N}) with \emph{N} the number of grid blocks. It should be noted that the size of the discretized contrast source \(\Delta\mitbf{\sigma E}\) needs to be padded by zeros such that the padded \(\Delta\mitbf{\sigma E}\) has twice the original number of points in all directions to avoid the periodicity in the FFT convolution result. The FFT of Green's tensor can be pre-calculated before calling the iterative solvers to save computational time during the iterative process.

\subsection{Domain Decomposition} \label{subsection_DD}

The domain decomposition method attempts to solve the problem for the entire from solutions of the different sub-domains \citep{saad_iterative_2003}. In our case, the spatial domain \(\Omega\) is decomposed into $M$~non-overlapping rectangular sub-domains \(\Omega_{j}\), hence
\begin{equation} \label{eq:DD_domains}
    \Omega = \bigcup_{j = 1}^{M}\Omega_{j},
\end{equation}
see Fig. \ref{fig:1_DD_schematic}. Adapting the domain decomposition formulation described in \citep{endo_large-scale_2009}, the convolution integral term or the scattered electric field term in equation (\ref{eq:E_IE}) can be expressed as a sum of scattered electric fields from each of the sub-domains. Subsequently, equation (\ref{eq:E_IE}) can be written as
\begin{equation} \label{eq:DD_E_IE}
    \mitbf{E}\left( \mitbf{r} \right) = \mitbf{E}^{(0)}\left( \mitbf{r} \right) + \sum_{j = 1}^{M}{\int_{\Omega_{j}}^{}{\mitbf{G}^{E}\ \left( \mitbf{r},\mitbf{r'} \right)\Delta\mitbf{\sigma}\left( \mitbf{r'} \right)\mitbf{E}\left( \mitbf{r'} \right) dV\left(\mitbf{r'}\right)}},
\end{equation}
where \(\Omega_{j}\) indicates the sub-domains with the conductivity anomaly. From equation (\ref{eq:DD_E_IE}), we obtain the following set of integral equations evaluated in each sub-domain:
\begin{equation} \label{eq:DD_discrete_E_IE}
    \mitbf{E}^{(i)} = \mitbf{E}^{(i,0)} + \sum_{j = 1}^{M}{\mitbf{\mathcal{G}}^{(ij)}\Delta\mitbf{\sigma}^{(j)}\mitbf{E}^{(j)}},\ \ i = 1,2,\ \ldots\ ,M.
\end{equation}
The terms \(\mitbf{E}^{(i,0)}\), \(\mitbf{E}^{(i)}\), and \({\Delta\mitbf{\sigma}}^{(i)}\) are the background electric field, total electric field, and the conductivity contrast defined at the sub-domain \(\Omega_{i}\), respectively. The terms \(\mathcal{G}^{(ij)}\Delta\mitbf{\sigma}^{(j)}\mitbf{E}^{(j)}\) in equation (\ref{eq:DD_discrete_E_IE}) are the discrete representations of the convolution integral in equation (\ref{eq:DD_E_IE}) which denote the scattered electric fields in the sub-domain \(\Omega_{i}\) due to the contrast source in the sub-domain \(\Omega_{j}\). It can be seen in equation (\ref{eq:DD_discrete_E_IE}) that the region without a conductivity anomaly does not contribute to the sum and hence can be omitted from the discretization when calculating the electric field. By collecting the scattered field terms into the left-hand side of the equations, the linear system of equations in (\ref{eq:DD_discrete_E_IE}) can be expressed with a block-matrix representation, viz.
\begin{equation} \label{eq:DD_AE_E}
    \mitbf{A}\widetilde{\mitbf{E}} = {\widetilde{\mitbf{E}}}^{(0)},
\end{equation}
where \(\mitbf{A}\) is the block matrix of the rearranged linear system according to the domain decomposition
\begin{equation} \label{eq:DD_A_definition}
    \mitbf{A} = 
    \begin{bmatrix}
\mitbf{I -}\mitbf{\mathcal{G}}^{(11)}\Delta\mitbf{\sigma}^{(1)} 
& - \mitbf{\mathcal{G}}^{(12)}\Delta\mitbf{\sigma}^{(2)}
& \ldots
& - \mitbf{\mathcal{G}}^{(1M)}\Delta\mitbf{\sigma}^{(M)} 
\cr 
- \mitbf{\mathcal{G}}^{(21)}\Delta\mitbf{\sigma}^{(1)}
& \mitbf{I -}\mitbf{\mathcal{G}}^{(22)}\Delta\mitbf{\sigma}^{(2)} 
& \ldots
& - \mitbf{\mathcal{G}}^{(2M)}\Delta\mitbf{\sigma}^{(M)} 
\cr 
\vdots
& \vdots
& \ddots
& \vdots
\cr 
- \mitbf{\mathcal{G}}^{(M1)}\Delta\mitbf{\sigma}^{(1)}
& - \mitbf{\mathcal{G}}^{(M2)}\Delta\mitbf{\sigma}^{(2)}
& \ldots
& \mitbf{I -}\mitbf{\mathcal{G}}^{(MM)}\Delta\mitbf{\sigma}^{(M)}  
\cr 
\end{bmatrix}\mathclose{},
\end{equation}
with \(\widetilde{\mitbf{E}}\) and \({\widetilde{\mitbf{E}}}^{(0)}\) are the block vectors containing the total and background electric fields in different sub-domains, respectively. These terms are defined as
\begin{equation} \label{eq:DD_E_E0_definition}
    \widetilde{\mitbf{E}} = 
    \begin{bmatrix}
\mitbf{E}^{(1)} \cr
\mitbf{E}^{(2)} \cr
 \vdots \cr
\mitbf{E}^{(M)} \cr
\end{bmatrix}\mathclose{},
\ \ \ \ \ {\widetilde{\mitbf{E}}}^{(0)} = \begin{bmatrix}
\mitbf{E}^{(1,0)} \cr
\mitbf{E}^{(2,0)} \cr
 \vdots \cr
\mitbf{E}^{(M,0)}
\end{bmatrix}\mathclose{}.
\end{equation}
Each block in the matrix \(\mitbf{A}\) indicates interaction terms between the sub-domains. The diagonal terms \(\left( \mitbf{I-}\mitbf{\mathcal{G}}^{(ii)}\Delta\mitbf{\sigma}^{(i)} \right)\) in equation (\ref{eq:DD_A_definition}) can be interpreted as the intra-domain interaction within a sub-domain while the off-diagonal terms \(-\mitbf{\mathcal{G}}^{(ij)}\Delta\mitbf{\sigma}^{(j)}\) represent the inter-domain interaction terms. Since the sub-domains are rectangular, the convolution integrals with Green's tensor in the intra- and inter-domain interaction terms can still be calculated using the FFT.

To solve the rearranged linear system of equation with domain decomposition in equation (\ref{eq:DD_AE_E}), the matrix \(\mitbf{A}\) is preconditioned by splitting the matrix into a strictly-lower-triangular \(\left( \mitbf{L} \right)\), strictly-upper-triangular \(\left( \mitbf{U} \right)\), and diagonal (\(\mitbf{D}\)) part \citep{barrett_templates_1994,saad_iterative_2003}:
\begin{equation} \label{eq:DD_A_split}
    \mitbf{A} = \left( \mitbf{L + U + D} \right),
\end{equation}
where the matrices \(\mitbf{L}\), \(\mitbf{U}\), and \(\mitbf{D}\) are defined by
\begin{equation} \label{eq:DD_D_L_U_definition}
\begin{aligned}   
\mitbf{L} &= \begin{bmatrix}
\mitbf{0} 
& \mitbf{0}
& \ldots
& \mitbf{0} 
\cr 
- \mitbf{\mathcal{G}}^{(21)}\Delta\mitbf{\sigma}^{(1)}
& \mitbf{0} 
& \ldots
& \mitbf{0} 
\cr 
\vdots
& \vdots
& \ddots
& \vdots
\cr 
- \mitbf{\mathcal{G}}^{(M1)}\Delta\mitbf{\sigma}^{(1)}
& - \mitbf{\mathcal{G}}^{(M2)}\Delta\mitbf{\sigma}^{(2)}
& \ldots
& \mitbf{0}  
\cr 
\end{bmatrix}\mathclose{}, \\ 
\mitbf{U} &= \begin{bmatrix}
\mitbf{0} 
& - \mitbf{\mathcal{G}}^{(12)}\Delta\mitbf{\sigma}^{(2)}
& \ldots
& - \mitbf{\mathcal{G}}^{(1M)}\Delta\mitbf{\sigma}^{(M)} 
\cr 
\mitbf{0}
& \mitbf{0} 
& \ldots
& - \mitbf{\mathcal{G}}^{(2M)}\Delta\mitbf{\sigma}^{(M)} 
\cr 
\vdots
& \vdots
& \ddots
& \vdots
\cr 
\mitbf{0}
& \mitbf{0}
& \ldots
& \mitbf{0}  
\cr 
\end{bmatrix}\mathclose{}, \\ 
\text{and} \\
\mitbf{D} &= \begin{bmatrix}
\mitbf{I -}\mitbf{\mathcal{G}}^{(11)}\Delta\mitbf{\sigma}^{(1)} 
& \mitbf{0}
& \ldots
& \mitbf{0} 
\cr 
\mitbf{0}
& \mitbf{I -}\mitbf{\mathcal{G}}^{(22)}\Delta\mitbf{\sigma}^{(2)} 
& \ldots
& \mitbf{0} 
\cr 
\vdots
& \vdots
& \ddots
& \vdots
\cr 
\mitbf{0}
& \mitbf{0}
& \ldots
& \mitbf{I -}\mitbf{\mathcal{G}}^{(MM)}\Delta\mitbf{\sigma}^{(M)}  
\cr 
\end{bmatrix}\mathclose{}, 
\end{aligned}
\end{equation}
respectively. By substituting the matrix splitting in (\ref{eq:DD_A_split}) into the equation (\ref{eq:DD_AE_E}) and some simple algebra, we obtain
\begin{equation} \label{eq:DD_original_form}
    {\widetilde{\mitbf{E}}} = \left( \mitbf{D} + \mitbf{L} \right)^{-1} \left\lbrack  {\widetilde{\mitbf{E}}}^{(0)}\mitbf{- U}{\widetilde{\mitbf{E}}}\right\rbrack\ \mathclose{},
\end{equation}
which can be solved by choosing an initial guess of ${\widetilde{\mitbf{E}}}$ and iteratively calculating the following
\begin{equation} \label{eq:DD_GS}
    {\widetilde{\mitbf{E}}}^{k + 1} = \left( \mitbf{D} + \mitbf{L} \right)^{- 1}\left\lbrack {\widetilde{\mitbf{E}}}^{(0)}\mitbf{- U}{\widetilde{\mitbf{E}}}^{k} \right\rbrack\mathclose{},
\end{equation}
with $k$ the iteration number. The iteration described in equation (\ref{eq:DD_GS}) corresponds to the block Gauss-Seidel iterative method \citep{barrett_templates_1994,saad_iterative_2003}. The matrix \(\left( \mitbf{D} + \mitbf{L} \right)\) has a lower triangular form where the inverse can be obtained using forward substitution \citep{venkateshan_chapter_2014}. The forward substitution process to compute equation (\ref{eq:DD_GS}) is outlined in Appendix \ref{appendix_forward_sub}. With the forward substitution, the total electric field update in each sub-domain according to equation (\ref{eq:DD_GS}) can be expressed in the simple form as:
\begin{equation} \label{eq:DD_GS_simpleform}
    \mitbf{E}^{(i),k + 1} = \left( \mitbf{I} - \mitbf{\mathcal{G}}^{(ii)}\Delta\mitbf{\sigma}^{(i)} \right)^{- 1}\left\lbrack \mitbf{E}^{(i,0)} + \sum_{j = 1}^{i - 1}{\mitbf{\mathcal{G}}^{(ij)}\Delta\mitbf{\sigma}^{(j)}\mitbf{E}^{(j),k + 1}} + \sum_{j = i + 1}^{M}{\mitbf{\mathcal{G}}^{(ij)}\Delta\mitbf{\sigma}^{(j)}\mitbf{E}^{(j),k}} \right\rbrack\mathclose{},
\end{equation}
where $i$ = 1,2,\ldots, \(M\) denotes the number of inner iterations where the IE is solved for one sub-domain and the number $k$ denotes the number of the total domain sweeps where the electric field is updated for the entire domain. The inverse operation of the block intra-domain term in equation (\ref{eq:DD_GS_simpleform}) is not calculated using the direct solver, but instead by using a Krylov subspace method to solve the following linear system of equations within each sub-domain:
\begin{equation} \label{eq:DD_GS_simple_Krylov}
    \left(\mitbf{I} - \mitbf{\mathcal{G}}^{(ii)}\Delta\mitbf{\sigma}^{(i)} \right)\mitbf{E}^{(i),k + 1} = \mitbf{E}^{(i,0)} + \sum_{j = 1}^{i - 1}{\mitbf{\mathcal{G}}^{(ij)}\Delta\mitbf{\sigma}^{(j)}\mitbf{E}^{(j),k + 1}} + \sum_{j = i + 1}^{M}{\mitbf{\mathcal{G}}^{(ij)}\Delta\mitbf{\sigma}^{(j)}\mitbf{E}^{(j),k}} \mathclose{}.
\end{equation}

The domain sweep is carried out until the relative residual on the whole domain reaches a desired threshold. The resulting operation of solving equation (\ref{eq:DD_GS_simple_Krylov}) iteratively is equivalent to the formulation described in \citet{zhdanov_integral_2006} and \citet{endo_large-scale_2009}. However, in our derivation, we can see the link between the original formulation to a block-preconditioned iterative method, which is the block Gauss-Seidel iterative method in this case. The convergence of the Gauss-Seidel iterative method depends on the diagonal dominance of the linear system matrix \citep{saad_iterative_2003}. In this case, if the sum of the inter-domain terms\textquotesingle{} norm is small compared to the norm of the intra-domain terms in equation (\ref{eq:DD_GS_simple_Krylov}), then this scheme is guaranteed to converge. Since the magnitude of Green's tensor elements depends on the distance between sub-domains, the interaction terms are small when the sub-domains are isolated from each other. When a sub-domain has small contrasts, the interaction is one-sided from the sub-domain with high contrast. These properties should be considered when designing the domain decomposition settings.

Instead of the Gauss-Seidel iterative method, one can also choose the Jacobi iterative method by taking only the diagonal part of the matrix \(\mitbf{A}\) as the preconditioner of the fixed-point equation instead of its lower triangular part. The fixed-point equation that corresponds to the Jacobi iterative method can be written as
\begin{equation} \label{eq:DD_Jacobi}
    {\widetilde{\mitbf{E}}}^{k + 1} = \mitbf{D}^{- 1}\left\lbrack {\widetilde{\mitbf{E}}}^{(0)} - \left( \mitbf{L}+\mitbf{U} \right){\widetilde{\mitbf{E}}}^{k} \right\rbrack\mathclose{},
\end{equation}
which leads to the following linear system of equations to be solved in each sub-domain:
\begin{equation} \label{eq:DD_Jacobi_Krylov}
    \left( \mitbf{I} - \mitbf{\mathcal{G}}^{(ii)}\Delta\mitbf{\sigma}^{(i)} \right)\mitbf{E}^{(i),k + 1} = \mitbf{E}^{(i,0)} + \sum_{j = 1, j \neq i}^{M}{\mitbf{\mathcal{G}}^{(ij)}\Delta\mitbf{\sigma}^{(j)}\mitbf{E}^{(j),k}} \mathclose{}.
\end{equation}
Since the right-hand side of equation (\ref{eq:DD_Jacobi_Krylov}) only depends on the solutions at the $k$-th iteration, the Jacobi iterative method is more straightforward to be implemented in parallel computing environments \citep{barrett_templates_1994}. In this case, the linear system of equations at each sub-domain can be solved with the Krylov solver in parallel and the interaction terms are updated after the Krylov solver computations are done in all sub-domains. The main drawback is that the Gauss-Seidel method generally has better convergence properties than the Jacobi method \citep{barrett_templates_1994}.

To further improve the computation speed, we propose to use a Krylov solver with adaptive target residual when solving the IE linear system of a sub-domain. The main idea is the relative residual of the Krylov solver in a sub-domain only needs to be an order of magnitude less than the full-domain relative residual to achieve the convergence of the Gauss-Seidel or Jacobi iteration. Inaccurate approximate solutions from the Krylov solver are acceptable at the beginning of the iteration and the relative residual target of the Krylov solver is lowered as the full-domain relative solver is decreasing during the Gauss-Seidel or Jacobi iteration. Additionally, the initial guess for the Krylov solver in the current outer iteration is updated from the result of the previous outer iteration. Detailed implementation of this strategy is shown in Algorithm 1.

To further accelerate the computation, it is possible to use the contraction integral form which accelerates the Krylov solver convergence rate \citep{endo_large-scale_2009,zhdanov_integral_2006}. However, in this study, we use the conventional integral equation formulation in the Krylov solver to reduce the complexity of evaluating the performance of the IE with domain decomposition.

\section{Numerical Results \& Discussion} \label{section_numerical_results}

In this section, we present three numerical cases to demonstrate the effectiveness of the domain decomposition preconditioning of the IE method. The first case is a model with two anomalous sub-domains separated by an isotropic medium with conductivity equal to the background conductivity. In the second case, we present a model where the anomalous isotropic conductivity is surrounded by an anisotropic medium. Lastly, we simulate a logging scenario across a faulted sand formation surrounded by anisotropic shale layers. We use the conventional IE formulation as described in section 2 in the GMRES solver for both full-domain IE and IE with domain decomposition (IE-DD) method. All numerical experiments presented in this paper are performed on a laptop with an AMD Ryzen 7 4800H processor and NVIDIA GeForce RTX 3060 Laptop GPU using MATLAB with GPU support enabled. We have compared our full-domain IE code with existing 1D semi-analytical solution \citep{shahriari2018numerical} and 3D finite-volume method \citep{hou2006finite}. This comparison is shown in Appendix \ref{appendix_verification} and our results show a good agreement with less than one per cent average absolute difference.

\subsection{Isolated Sub-domains Example} \label{subsection_isolated_subdom_ex}

We consider two isolated anomalous sub-domains embedded in an isotropic medium background as shown in Fig. \ref{fig:2_model_example1}. The background conductivity \(\sigma_{0}\) is equal to 0.1~S~m$^{-1}$ and the conductivity in the anomalous sub-domain is equal to 0.01~S~m$^{-1}$. A transmitter with 24 KHz frequency is located at the origin (x = 0, y = 0, z = 0~m) and it is oriented in the x-direction. The whole domain is discretized into 128 $\times$ 128 $\times$ 128 grid blocks with a uniform grid size of 0.25 $\times$ 0.25 $\times$ 0.25~m\textsuperscript{3}. The two anomalous sub-domains are set to have an equal size of 30 $\times$ 30 $\times$ 7.5~m\textsuperscript{3}. The distance between the closest edges of the two sub-domains is 10~m which is approximately equal to the skin depth of the background medium given the transmitter frequency.

To solve the conventional full-domain IE and the inverse operation in IE-DD formulation with Gauss-Seidel iterative method, we use restarted GMRES method with 10 restart iterations. We set the relative residual $e$ = 10\textsuperscript{-6} for both the conventional full-domain IE and IE-DD Gauss-Seidel iteration stopping criterion. For this case, we did not implement the adaptive relative residual scheme and set $e$ = 10\textsuperscript{-6} as the relative residual target for the GMRES solver stopping criterion in the IE-DD to analyze the convergence behavior of the method. Because the medium without contrast does not contribute to the scattering field, the relative residuals are only computed within the anomalous sub-domains in both cases. 

The sub-domains without the conductivity anomalies are excluded from the discretization in the IE-DD iterations. This results in only half of the total number of grid blocks of the full-domain IE being discretized in the IE-DD iterations. Fig.~\ref{fig:3_ex1_relres_iter} shows the convergence plot and the total number of GMRES iterations taken to reach the target residual. The number of GMRES iterations required is decreasing in each of the Gauss-Seidel iterations as the relative residual is decreasing. This indicates that the changes in the electric fields due to the interaction terms become smaller as the initial guess for the GMRES solver is updated in each of the Gauss-Seidel iterations. Overall, it took only four Gauss-Seidel iterations for the method to converge below the threshold level with a total of 278 GMRES iterations and the total computational time is 23~s. The full-domain IE took 86 GMRES iterations to converge on the same error level, but the computation time is 53~s. Even though there is more GMRES iteration in IE-DD, the total computational time is approximately twice faster compared to the full-domain IE because the operation at each of the GMRES iterations in IE-DD works on a smaller domain with the number of blocks equal to a quarter of the full-domain grid blocks in each domain.

The magnetic field comparison between both methods is shown in Fig.~\ref{fig:4_ex1_fields}. Qualitatively, there are no differences observed because both methods show similar numerical results within less than 0.01~per~cent average normalized magnitude difference. Therefore, the IE-DD method will give the same result within the same relative residual level as the full-domain IE.

\subsection{Simple Anisotropic Medium Example} \label{subsection_simple_anisotropic_ex}

Fig.~\ref{fig:5_ex2_model} shows an xz-plane view of a faulted resistive isotropic medium with a conductivity of 0.01~S~m$^{-1}$ surrounded by an anisotropic medium with vertical transverse isotropy. The conductivity tensor of the anisotropic medium consists of the conductivity in the horizontal and vertical direction with the value of \(\sigma_{h}\) = 0.2~S~m$^{-1}$ and \(\sigma_{v}\) = 0.1~S~m$^{-1}$, respectively. The conductivity of the media does not vary in the y-direction. For the background medium, we choose a homogenous isotropic medium with the conductivity of \(\sigma_{0}\) = 0.01~S~m$^{-1}$. A transmitter with 24 KHz frequency is located at the origin and is oriented in the x-direction. We set 10\textsuperscript{-6} as the relative residual target for both methods. The whole domain is discretized into 120 $\times$ 120 $\times$ 120 grid blocks with a grid size of 0.25 $\times$ 0.25 $\times$ 0.25~m\textsuperscript{3}.

The full domain is decomposed into three rectangular sub-domains of equal size as illustrated in Fig.~\ref{fig:5_ex2_model}. Each sub-domain is discretized into 120 $\times$ 120 $\times$ 40~grid blocks with a grid size of 0.25 $\times$ 0.25 $\times$ 0.25~m\textsuperscript{3}. With this decomposition, the faulted resistive layer is located only in sub-domain 1 while the other two sub-domains contain only the anisotropic medium.

We present three different schemes of IE-DD to calculate the electric field of the model. The first one (IE-DD-GS-F) is the IE-DD with Gauss-Seidel iterative method and fixed GMRES solver relative residual stopping criterion equal to 10\textsuperscript{-6} in every Gauss-Seidel iteration. The second one (IE-DD-GS-A) is the IE-DD with Gauss-Seidel iterative method and adaptive GMRES solver relative residual stopping criterion. The third one (IE-DD-Jacobi-A) is the IE-DD with Jacobi iterative method and with the same adaptive GMRES relative residual as the second scheme. In the adaptive relative residual scheme, the relative residual stopping criterion is set to be one order of magnitude lower than the relative residual calculated on the whole domain or full-domain relative residual divided by ten.

Fig.~\ref{fig:6_ex2_total_GMRES_iter} displays the comparison between the total GMRES iteration in each Gauss-Seidel iteration for both IE-DD-GS schemes. The total number of GMRES iterations generally increases along with the Gauss-Seidel iteration in the adaptive relative residual scheme, while it is decreasing in the fixed relative residual scheme. In both cases, sub-domain 1, which contains the faulted resistive layer, took the greatest number of GMRES iterations. Since the number of GMRES iterations is proportional to the conductivity contrast, this indicates that sub-domain 1 has the largest conductivity contrast compared to the other two sub-domains. 

Table~\ref{table:schemes_performance} summarizes the computational cost comparison between the full domain IE and IE-DD with three different schemes. Based on the computation time, the IE-DD-GS-A is the fastest scheme compared to the other two methods. The IE-DD-GS-A scheme shows a faster computation time compared to the IE-DD-GS-F because there are fewer GMRES iterations in the IE-DD-GS-A scheme. Therefore, specifying the adaptive relative residual for the Krylov solver in the IE-DD improves the computation time of the original IE-DD formulation with the cost of going through more Gauss-Seidel iterations. The full-domain relative residual plot in each of the outer iterations shown in Fig.~\ref{fig:7_ex2_fullrelres_GS_Jacobi} indicates that the IE-DD-GS-A has a better convergence rate compared to the IE-DD-Jacobi-A. Besides having a smaller number of total outer iterations, the number of GMRES iterations taken in the IE-DD-GS-A is also fewer compared to the IE-DD-Jacobi-A. However, the computation time of the IE-DD-Jacobi-A can potentially be improved by further utilizing parallel computation to independently solve the linear system in each sub-domains.

\subsection{Logging Simulation Across a Complex Formation} \label{subsection_logging_ex}

We simulated induction logs across the faulted anisotropic formation with an 85$^\circ$ drilling angle as illustrated in Fig.~\ref{fig:8_ex3_model}a. This formation consists of anisotropic shale layers surrounding isotropic sand layers. The shale layers are indicated by the blue colours and the sand layers are indicated by the yellow colours in Fig. \ref{fig:8_ex3_model}. The model has 2.5D main structural features with the addition of a simple 3D Gaussian perturbation only in the sand layers to imitate a fluid distribution in a reservoir. This perturbation is defined by
\begin{equation} \label{eq:Gaussian_perturb}
    \mitbf{\sigma}_{sand} = \mitbf{ \sigma}_{sand}^{u} + \alpha\mitbf{ \sigma}_{sand}^{u}\exp\left( - \frac{\left| \mitbf{r}_{sand}-\mitbf{r}_{c} \right|}{\gamma} \right)\mathclose{},
\end{equation}
where the subscripts $sand$ denote the values located in the sand layers and the superscripts $u$ indicate the defined unperturbed value; \(\mitbf{r}_{c}\) is the location of the maximum perturbation; \(\alpha\) and \(\gamma\) are the factors that control the magnitude and range of the perturbation, respectively. In this example, we set the peak perturbation location \(\mitbf{r}_{c}\) at x = 500~m, y = 0~m, and z = 40~m; and define \(\alpha\) = 4 and \(\gamma\) = 50~m.

We use a moving 3D forward modelling window to simulate a moving transmitter scenario. The z-direction in the forward modelling window is directed to the drilling direction so it is consistent with the component direction of the induction tools \citep{pardo_fast_2015}. Hence the coordinate system in the window is rotated from the cartesian coordinate according to the drilling direction as illustrated in Fig.~\ref{fig:8_ex3_model}a. Consequently, the conductivity tensor elements are transformed following the domain rotation \citep{gao_simulation_2006}, see Appendix \ref{appendix_conductivity} for detail. In each of the forward modelling windows, we set a constant background conductivity \(\sigma_{0}\) = 0.1~S~m$^{-1}$.

Following the typical tool configurations described in \citet{antonsen_what_2022}, we set a z-oriented transmitter with a frequency of 24 KHz and three receivers with spacings of 7, 15, and 30~m as illustrated in Fig.~\ref{fig:9_tools_schematic} for the logging simulations. A forward modelling window with a size of 32 $\times$ 32 $\times$ 32~m\textsuperscript{3} may not be enough to capture the full sensitivity of all the receivers, especially the one with the largest spacing. Hence, we tested two different windows with different sizes of 32 $\times$ 32 $\times$ 64~m\textsuperscript{3} and 64 $\times$ 64 $\times$ 64~m\textsuperscript{3} to see different sensitivities of the receivers with the forward modelling domain size. We refer to the smaller window as window 1 and the larger one as window 2. In both windows, we keep a grid size of 0.25 $\times$ 0.25 $\times$ 0.25~m\textsuperscript{3} resulting in a total of 128 $\times$ 128 $\times$ 256 and 256 $\times$ 256 $\times$ 256~grid blocks for window 1 and 2, respectively. Since our current computer GPU memory size can only handle a maximum of 128\textsuperscript{3}~grid blocks for the GMRES solver computation, we decompose the window 1 and 2 into two and eight sub-domains in the drilling direction, respectively, as illustrated in Fig.~\ref{fig:10_windows_sketch}. In this way, the memory requirements to calculate the electric field with both windows are reduced to the memory requirement for the calculation using 128\textsuperscript{3}~grid blocks plus the memory for storing the electric fields. This allows us to fully take advantage of the acceleration with the GPU implementation. 

The logging position starts at x = 0~m, y = 0~m, z = 0~m and ends at x = 900~m, y = 0~m, z = 78.74~m. In each logging position, we use the IE-DD-GS-A scheme and set a full-domain tolerance of 10\textsuperscript{-3} as the stopping criterion. Fig.~\ref{fig:11_ex3_log_response} shows the magnitude of the z-component of the magnetic field \textbar{}\(H_{zz}\)\textbar{} measured in the receivers at each transmitter position. Qualitatively, the differences in the results between the two window settings are increasing with the receiver spacings. This result shows different sensitivities of the transmitter-receiver configuration, and we can observe that the sensitivity range is proportional to the receiver spacing.

The computation time required to calculate the magnetic field for one logging position using the window 1 and 2 settings takes an average of approximately 1.5 and 15~minutes, respectively. Updating the interaction terms is the most expensive part of the computation time, taking up around 80 per cent of the time at every iteration due to the operation acting on the entire domain that consists of a massive number of grid blocks. In every position, it took less than seven Gauss-Seidel iterations to reach the desired tolerance. 

\section{Conclusion} \label{section_conclusion}

The linear system of equations arising in the IE method for 3D EM method modelling can be naturally decomposed into a set of linear systems of equations that correspond to the IE in different parts of the modelling domain. The IE-DD formulation is reducing the memory requirement to compute a large-scale problem as it provides the connection between each sub-domains while still maintaining the viability of using FFTs to calculate the convolution integral operation. By expressing these linear systems of equations in a block matrix representation where each block represents the interactions between the domains, we have made a link between the derivation in \citet{zhdanov_integral_2006} and \citet{endo_large-scale_2009} with a preconditioned fixed-point iteration using domain decomposition method. Depending on the choice of the preconditioner, the fixed-point iteration corresponds to the block Gauss-Seidel and Jacobi iterative method. In every Gauss-Seidel or Jacobi iteration, the inverse of the block intra-domain interaction term is calculated using the Krylov subspace method instead of a direct solver. 

Our numerical experiment results show that a reduction in computation time can be achieved although the total number of GMRES solver iterations in IE-DD schemes is more than in the full-domain IE. This speed-up is due to the GMRES solver in the decomposed domains being cheaper to compute and it is shown that it only takes five to eight IE-DD outer iterations to reach the desired tolerance. Additionally, specifying adaptive relative residual stopping improves the computation time of the IE-DD by reducing the total number of GMRES iterations required for reaching desired error tolerance. The Gauss-Seidel preconditioning with adaptive relative residual has the fastest computation time among the schemes that we tested in this study. This scheme reduces the computation time of the conventional IE by approximately 35 per cent. The scheme with Jacobi preconditioning takes longer computation time compared to the one with Gauss-Seidel. However, the form of the Jacobi iterative method is more suitable for parallel computation as the operation in each sub-domain can be computed independently, which is a subject for future implementation.

In this study, we have only implemented IE-DD with a simple iterative update corresponding to Gauss-Seidel and Jacobi iterative method. The Gauss-Seidel and Jacobi iterative methods are in general not very competitive in terms of convergence compared to the Krylov subspace method \citep{barrett_templates_1994}. Therefore, further potential improvement of the IE-DD presented in this study is obtained by implementing the Krylov subspace as the outer iteration update instead of the Gauss-Seidel and Jacobi iteration update. Another interesting application of the domain decomposition in the IE method would be to incorporate a direct method that can be computed in parallel into the domain decomposition preconditioner, for example using the T-matrix method \citep{jakobsen_distorted_2018,sommer2018towards}.

\section{Acknowledgements}

This work is part of the Center for Research-based Innovation DigiWells: Digital Well Center for Value Creation, Competitiveness and Minimum Environmental Footprint (NFR SFI project no. 309589, https://DigiWells.no). The center is a cooperation of NORCE Norwegian Research Centre AS, the University of Stavanger, the Norwegian University of Science and Technology (NTNU), and the University of Bergen. It is funded by Aker BP, ConocoPhillips, Equinor, Lundin Energy, TotalEnergies, Vår Energi, Wintershall Dea, and the Research Council of Norway.

\section{Data Availability}
Currently, the data and code relating to this work are not freely available. We are considering publishing the codes with an open-source license in the future.

\bibliography{references_gji}

\appendix
\section{Forward-Substitution in Gauss-Seidel Iterative
Method} \label{appendix_forward_sub}

The inverse of a lower-triangular matrix can be obtained through forward substitution (Venkateshan \& Swaminathan, 2014). For simplicity, we demonstrate the forward-substitution process in the case of domain decomposition with three sub-domains. We can write the block-matrix representation and its splitting for three sub-domains:
\begin{equation} \label{eq:Appendix_A3domains}
    A = \begin{bmatrix}
\mitbf{A}_{11} & \mitbf{A}_{12} & \mitbf{A}_{13} \\
\mitbf{A}_{21} & \mitbf{A}_{22} & \mitbf{A}_{23} \\
\mitbf{A}_{31} & \mitbf{A}_{32} & \mitbf{A}_{33} \\
\end{bmatrix} = \left( \mitbf{D + L + U} \right),
\end{equation}
with the matrices \mitbf{D}, \mitbf{L}, and \mitbf{U} are the block diagonal, strictly lower-triangular, and strictly upper-triangular of \mitbf{A}, respectively. These terms are defined as
\begin{equation} \label{eq:Appendix_DLU}
\begin{aligned}    
\mitbf{D} &= \begin{bmatrix}
\mitbf{A}_{11} & \mitbf{0} & \mitbf{0} \\
\mitbf{0} & \mitbf{A}_{22} & \mitbf{0} \\
\mitbf{0} & \mitbf{0} & \mitbf{A}_{33} \\
\end{bmatrix}\mathclose{}, \\
\mitbf{L} &= \begin{bmatrix}
\mitbf{0} & \mitbf{0} & \mitbf{0} \\
\mitbf{A}_{21} & \mitbf{0} & \mitbf{0} \\
\mitbf{A}_{31} & \mitbf{A}_{32} & \mitbf{0} \\
\end{bmatrix}\mathclose{}, \\
\text{and}\\
\mitbf{U} &= \begin{bmatrix}
\mitbf{0} & \mitbf{A}_{12} & \mitbf{A}_{13} \\
\mitbf{0} & \mitbf{0} & \mitbf{A}_{23} \\
\mitbf{0} & \mitbf{0} & \mitbf{0} \\
\end{bmatrix}\mathclose{}.
\end{aligned}
\end{equation}
Here each of the blocks denotes the inter-domain and intra-domain operators described in section~\ref{section_theory} as
\begin{gather}
    \mitbf{A}_{ii} = \left( \mitbf{I -}\mitbf{\mathcal{G}}^{(ii)}\Delta\mitbf{\sigma}^{(i)} \right),\label{eq:Appendix_intradomain}\\
    \mitbf{A}_{ij} = -\mitbf{\mathcal{G}}^{(ij)}\Delta\mitbf{\sigma}^{(j)},\ \ i \neq j \label{eq:Appendix_interdomain}
\end{gather}
The fixed-point equations using matrix \textbf{A} that corresponds to the Gauss-Seidel iterative method
\begin{equation} \label{eq:Appendix_GS_generalform}
    {\widetilde{\mitbf{E}}}^{k + 1} = \left( \mitbf{D} + \mitbf{L} \right)^{-1}\left\lbrack \mitbf{E}^{(0)}\mitbf{-U}{\widetilde{\mitbf{E}}}^{k} \right\rbrack.
\end{equation}
By substituting the matrices \mitbf{D}, \mitbf{L}, and
\mitbf{U} for 3 sub-domains and calculating the inverse of
\(\left( \mitbf{D} + \mitbf{L} \right)\) , we obtain:
\begin{equation} \label{eq:Appendix_GS_matrix}
    \begin{bmatrix}
\mitbf{E}^{(1)} \\
\mitbf{E}^{(2)} \\
\mitbf{E}^{(3)} \\
\end{bmatrix}^{k + 1} = \begin{bmatrix}
\mitbf{A}_{11}^{-1} & \mitbf{0} & \mitbf{0} \\
 - \mitbf{A}_{22}^{-1}\mitbf{A}_{21}\mitbf{A}_{11}^{-1} & \mitbf{A}_{22}^{-1} & \mitbf{0} \\
\mitbf{- A}_{31}^{-1} & - \mitbf{A}_{33}^{-1}\mitbf{A}_{32}\mitbf{A}_{22}^{-1} & \mitbf{A}_{33}^{-1} \\
\end{bmatrix}\begin{bmatrix}
\mitbf{R}_{1} \\
\mitbf{R}_{2} \\
\mitbf{R}_{3} \\
\end{bmatrix}\mathclose{},
\end{equation}
where
\begin{equation}
    \mitbf{A}_{31}^{-1} = \mitbf{A}_{33}^{-1}\mitbf{A}_{31}\mitbf{A}_{11}^{-1} -\mitbf{A}_{33}^{-1}\mitbf{A}_{32}\mitbf{A}_{22}^{-1}\mitbf{A}_{21}\mitbf{A}_{11}^{-1},
\end{equation}
and the terms \(\mitbf{R}_{i}\) denote the $i$-th row of the second term in the right-hand side of equation (\ref{eq:Appendix_GS_generalform}) written as
\begin{equation}
    \begin{bmatrix}
\mitbf{R}_{1} \\
\mitbf{R}_{2} \\
\mitbf{R}_{3} \\
\end{bmatrix} = \begin{bmatrix}
\mitbf{E}^{(1,0)} - \mitbf{A}_{12}\mitbf{E}^{(2),k} - \mitbf{A}_{13}\mitbf{E}^{(3),k} \\
\mitbf{E}^{(2,0)} - \mitbf{A}_{23}\mitbf{E}^{(3),k} \\
\mitbf{E}^{(3,0)} \\
\end{bmatrix}\mathclose{}.
\end{equation}
By multiplying the matrix on the right-hand side of equation (A.6), we obtain the following equations:
\begin{gather}
    \mitbf{E}^{(1),k + 1} = \mitbf{A}_{11}^{- 1}\mitbf{R}_{1},\label{eq:Appendix_GS_dom1}\\
    \mitbf{E}^{(2),k + 1} = \mitbf{A}_{22}^{- 1}\left\lbrack \mitbf{R}_{2} - \mitbf{A}_{21}\mitbf{A}_{11}^{-1}\mitbf{R}_{1} \right\rbrack\mathclose{}, \label{eq:Appendix_GS_dom2}\\
    \mitbf{E}^{(3),k + 1} = \mitbf{A}_{33}^{- 1}\left\lbrack \mitbf{R}_{3} - \mitbf{A}_{32}\mitbf{A}_{22}^{-1}\left( \mitbf{R}_{2} - \mitbf{A}_{21}\mitbf{A}_{11}^{-1}\mitbf{R}_{1} \right) - \mitbf{A}_{31}\mitbf{A}_{11}^{- 1}\mitbf{R}_{1} \right\rbrack\mathclose{}.\label{eq:Appendix_GS_dom3}
\end{gather}

Notice that the term \(\mitbf{A}_{11}^{-1}\mitbf{R}_{1}\) in the second right-hand side term of equations (\ref{eq:Appendix_GS_dom2}) and the third right-hand side term of (\ref{eq:Appendix_GS_dom3}) can be substituted by \(\mitbf{E}^{(1),k + 1}\) from equation (\ref{eq:Appendix_GS_dom1}). Also, the term \(\mitbf{A}_{22}^{-1}\left( \mitbf{R}_{2}\ -\mitbf{A}_{21}\mitbf{A}_{11}^{-1}\mitbf{R}_{1} \right)\) in the second right-hand side term of equation (\ref{eq:Appendix_GS_dom3}) can be substituted by \(\mitbf{E}^{(2),k + 1}\) from equation (\ref{eq:Appendix_GS_dom2}). With substitutions on these terms, equations (\ref{eq:Appendix_GS_dom1}-\ref{eq:Appendix_GS_dom3}) can be expressed as:
\begin{gather}
    \mitbf{E}^{(1),k + 1} = \mitbf{A}_{11}^{-1}\mitbf{R}_{1}, \label{eq:Appendix_GSsubs_dom1}\\
    \mitbf{E}^{(2),k + 1} = \mitbf{A}_{22}^{-1}\left\lbrack \mitbf{R}_{2} - \mitbf{A}_{21}\mitbf{E}^{(1),k + 1} \right\rbrack\mathclose{},\label{eq:Appendix_GSsubs_dom2}\\
    \mitbf{E}^{(3),k + 1} = \mitbf{A}_{33}^{-1}\left\lbrack \mitbf{R}_{3} -\mitbf{A}_{32}\mitbf{E}^{(2),k + 1} - \mitbf{A}_{31}\mitbf{E}^{(1),k + 1} \right\rbrack\mathclose{},\label{eq:Appendix_GSsubs_dom3}
\end{gather}
where the terms \(\mitbf{A}_{ii}^{-1}\) are the block matrices that indicate the process of solving a linear system of equations in the $i$-th sub-domain. It can be observed from the equations (\ref{eq:Appendix_GSsubs_dom1}-\ref{eq:Appendix_GSsubs_dom3}) that, in general, there are recurrences of the term \(\mitbf{E}^{(i),k + 1}\) in all of the equations in the $j$-th sub-domain for $j$\textgreater $i$. This implies the results of the fixed-point equation (\ref{eq:Appendix_GS_generalform}) can be obtained by sequentially solving the linear system of equations in each sub-domain and updating the right-hand side in the equations for the next sub-domain using the most recent solutions.

\section{Verification of the Conventional IE-Code} \label{appendix_verification}

To verify the accuracy of our conventional 3D IE code, we compare our numerical results to those obtained with a semi-analytical 1D solver \citep{shahriari2018numerical} and a 3D finite-volume solver \citep{hou2006finite}. Following \citep{nazanin_limitdetect}, we consider a logging-while-drilling simulation across a layered anisotropic medium with a drilling angle of 80$^\circ$ as illustrated in Fig.~\ref{fig:appendix1_model}. The tool consists of a tri-axial transmitter and receiver with a frequency of 12~KHz and a receiver spacing of 7.62~m. We use a moving forward modelling window with a size of 48.64 $\times$ 48.64 $\times$ 48.64~m$^3$. This window is discretized into 128 $\times$ 128 $\times$ 128 grid blocks and a cell size of 0.38 $\times$ 0.38 $\times$ 0.38~m$^{3}$. We set a constant background conductivity of 0.1118 S m$^{-1}$.

Fig.~\ref{fig:appendix2_comparison} shows the co-axial and co-planar components of the magnetic fields obtained using different numerical methods. We do not include the Y co-planar components of the magnetic fields because the values are zero. Overall, the results obtained from our IE-code show a good agreement with the 1D semi-analytical and 3D finite volume results. The average absolute difference of all components calculated using our IE code is less than one per cent compared to the 1D semi-analytical result.

\section{Conductivity Tensor Transformation for Transversely
Anisotropic Formation} \label{appendix_conductivity}

The tensor structure of electrical conductivity \(\mitbf{\sigma}\) for
anisotropic media is generally expressed as \citep{zhdanov_geophysical_2009}
\begin{equation}
    \mitbf{\sigma} = \begin{bmatrix}
\sigma_{xx} & \sigma_{xy} & \sigma_{xz} \\
\sigma_{yx} & \sigma_{yy} & \sigma_{yz} \\
\sigma_{zx} & \sigma_{zy} & \sigma_{zz} \\
\end{bmatrix}\mathclose{},
\end{equation}
with an off-diagonal symmetry \(\sigma_{ij} = \sigma_{ji}\). For a
vertical transverse isotropic medium with the z-axis as the vertical axis, the conductivity tensor is written as \citep{gao_simulation_2006,jakobsen_distorted_2018}\begin{equation}
\mitbf{\sigma} = \begin{bmatrix}
\sigma_{h} & 0 & 0 \\
0 & \sigma_{h} & 0 \\
0 & 0 & \sigma_{v} \\
\end{bmatrix}\mathclose{},
\end{equation}
where \(\sigma_{h}\) and \(\sigma_{v}\) are the conductivity in the horizontal and vertical directions, respectively. In the case where the angle between the formation layering and the drilling trajectory is not
90, it is necessary to rotate the conductivity tensor from the formation coordinate system to the induction tool coordinate system. For a transversely isotropic formation, the explicit expressions of the rotated conductivity tensor are as follows \citep{gao_simulation_2006}
\begin{equation}
    \begin{aligned}
        \sigma_{xx}' &= \sigma_{h} + \left( \sigma_{v} - \sigma_{h} \right)\sin^{2}\theta\cos^{2}\phi,
\\
\sigma_{xy}' &= \left( \sigma_{v} - \sigma_{h} \right)\sin^{2}\theta\sin\phi\cos\phi,
\\
\sigma_{xz}' &= \left( \sigma_{v} - \sigma_{h} \right)\sin\theta\cos\theta\cos\phi,
\\
\sigma_{yy}' &=\sigma_{h} + \left( \sigma_{v} - \sigma_{h} \right)\sin^{2}\theta\sin^{2}\phi,
\\
\sigma_{yz}' &=\left( \sigma_{v} - \sigma_{h} \right)\sin\theta\cos\theta\sin{\phi},
\\ \text{and} \\
\sigma_{zz}' &=\sigma_{v} - \left( \sigma_{v} - \sigma_{h} \right)\sin^{2}\theta,
    \end{aligned}
\end{equation}
where \(\theta\) and \(\phi\) are the z-axis rotation and y-axis rotation angles from the formation coordinate system to the induction tool coordinate system, respectively. The resulting conductivity tensor has non-zero off-diagonal components and our method can deal with this complication without making any other changes in the implementation.

\newpage
\FloatBarrier
\section*{Figures, Table, Algorithm}

\renewcommand{\thefigure}{\arabic{figure}}

\begin{figure}
    \centering
    \includegraphics[,width=7cm,height=5cm]{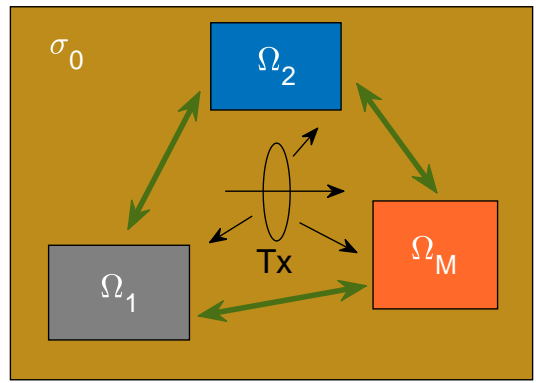}
    \caption{Schematic of the domain decomposition. The black arrows indicate the field coming from the transmitter Tx to the sub-domains \(\Omega_{j}\). The double-headed green arrows indicate the scatterers\textquotesingle{} interaction between the sub-domains. The transmitter can also be located in the anomalous domain.}
    \label{fig:1_DD_schematic}
\end{figure}

\begin{figure}
    \centering
    \includegraphics[width=15cm,height=5.6cm]{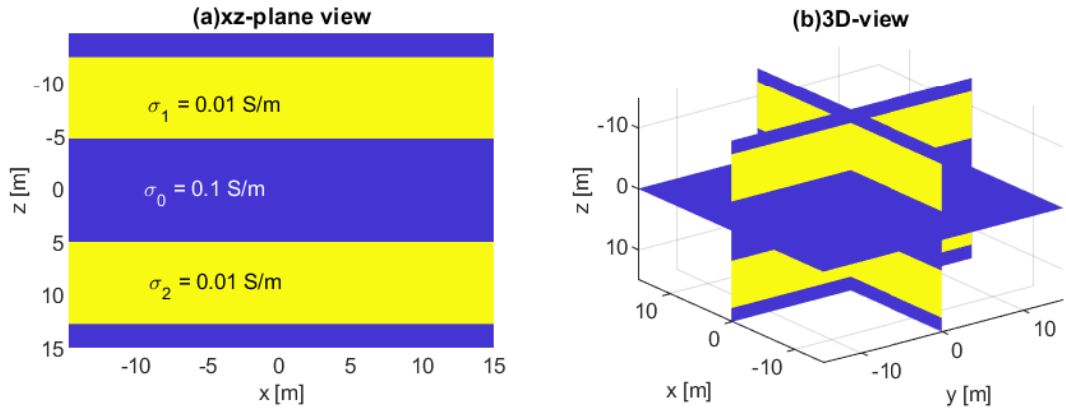}
    \caption{(a) xz-plane view of the model at y = 0 m. (b) xy- and xz-slice of the model in 3D view.}
    \label{fig:2_model_example1}
\end{figure}

\begin{figure}
    \centering
    \includegraphics[width=15cm,height=5.6cm]{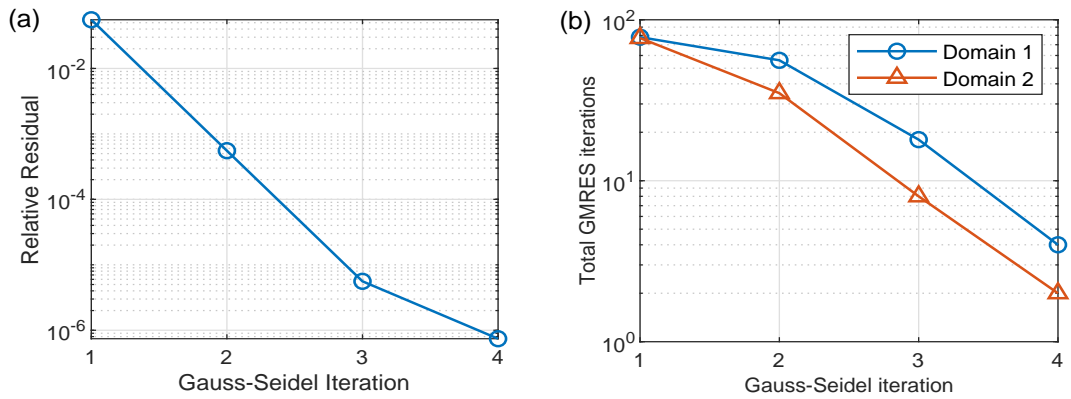}
    \caption{(a) The relative residual including the anomalous domains interaction and (b) the number of GMRES solver iterations within each of the Gauss-Seidel iterations.}
    \label{fig:3_ex1_relres_iter}
\end{figure}

\begin{figure}
    \centering
    \includegraphics[width=15cm,height=11cm]{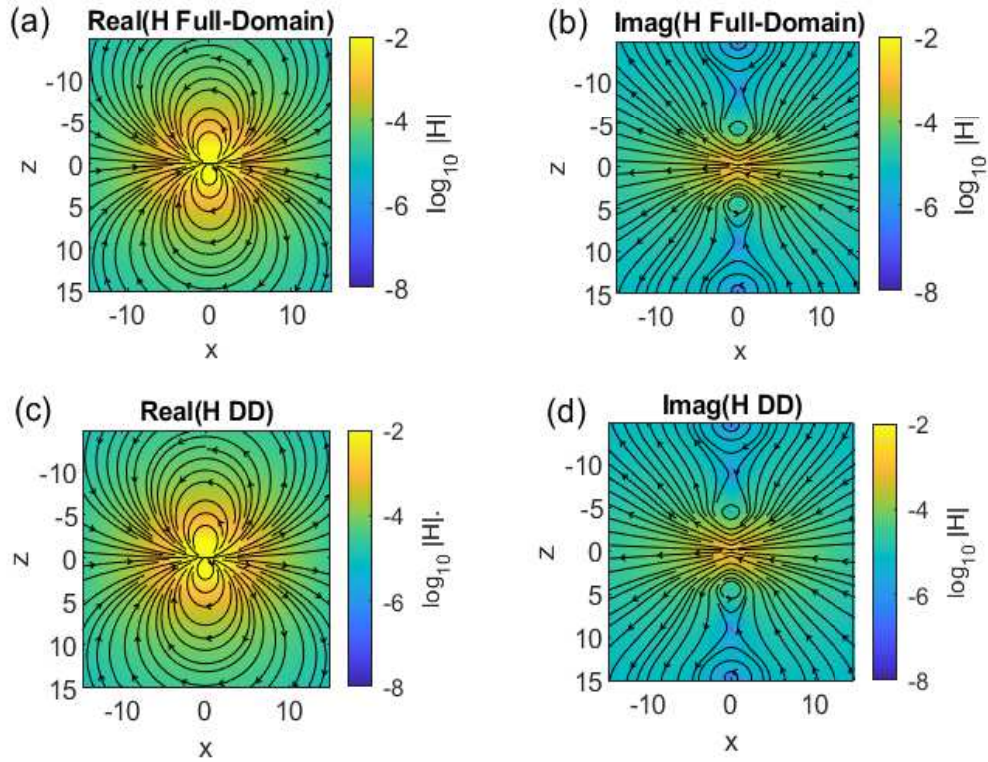}
    \caption{xz-plane slice of magnetic fields at y = 0~m. (a) Real and (b) imaginary parts of the magnetic field obtained from the full-domain IE. (c) Real and (d) imaginary parts of the magnetic field obtained from the IE-DD. A transmitter oriented to the x-direction is located at x = 0~m, y = 0~m, and z = 0~m.}
    \label{fig:4_ex1_fields}
\end{figure}

\begin{figure}
    \centering
    \includegraphics[width=15cm,height=6cm]{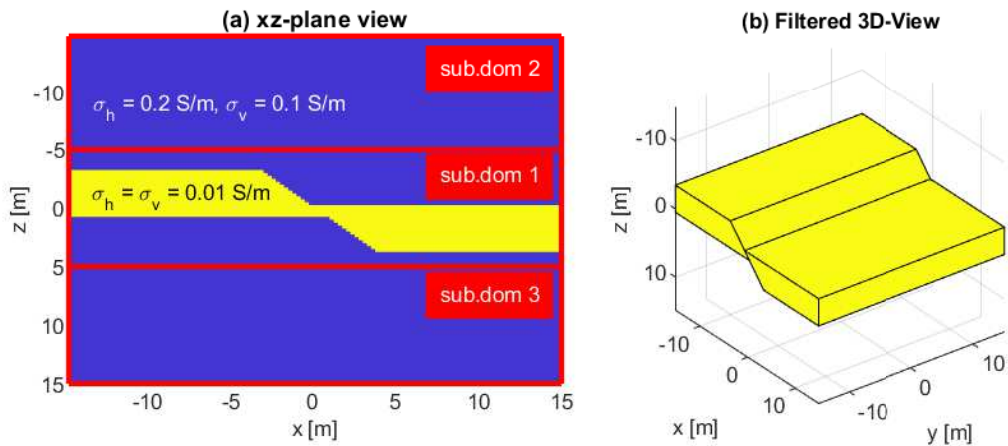}
    \caption{(a) xz-plane view of a faulted resistive layer surrounded by an anisotropic medium at y = 0~m and the domain decomposition setting. (b) 3D view with the anisotropic layers removed. A transmitter oriented to the x-direction is located at x = 0~m, y = 0~m, and z = 0~m.}
    \label{fig:5_ex2_model}
\end{figure}

\begin{figure}
    \centering
    \includegraphics[width=15cm,height=7.5cm]{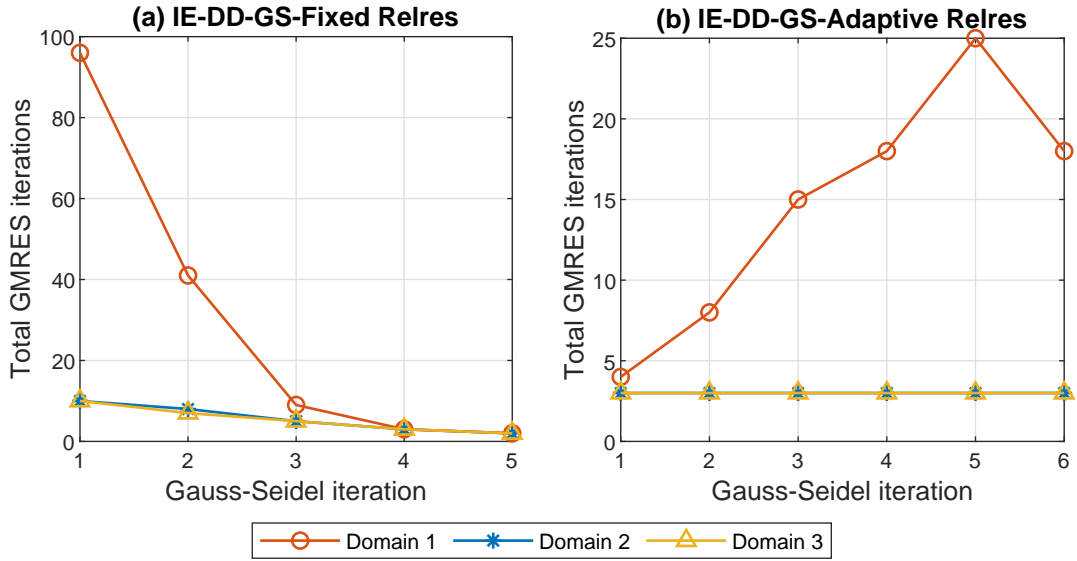}
    \caption{ Total GMRES iterations within each Gauss-Seidel iteration for IE-DD-GS with (a) fixed and (b) adaptive relative residual stopping criterion.}
    \label{fig:6_ex2_total_GMRES_iter}
\end{figure}

\begin{figure}
    \centering
    \includegraphics[width=7.5cm,height=5.5cm]{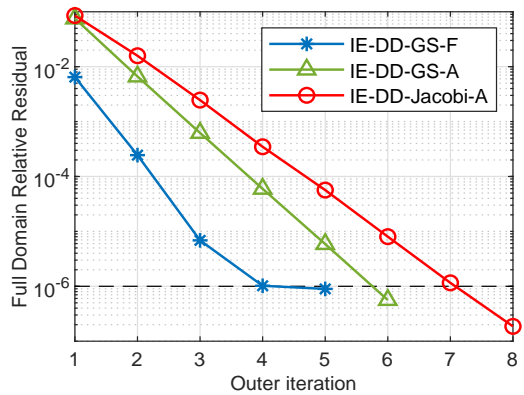}
    \caption{The convergence rate comparison between IE-DD-GS-F (blue line), IE-DD-GS-A (green line), and IE-DD-Jacobi-A (red line). The black dashed line indicates the stopping criterion for the Gauss-Seidel and Jacobi iterations.}
    \label{fig:7_ex2_fullrelres_GS_Jacobi}
\end{figure}

\begin{figure}
    \centering
    \includegraphics[width=14.67cm,height=11cm]{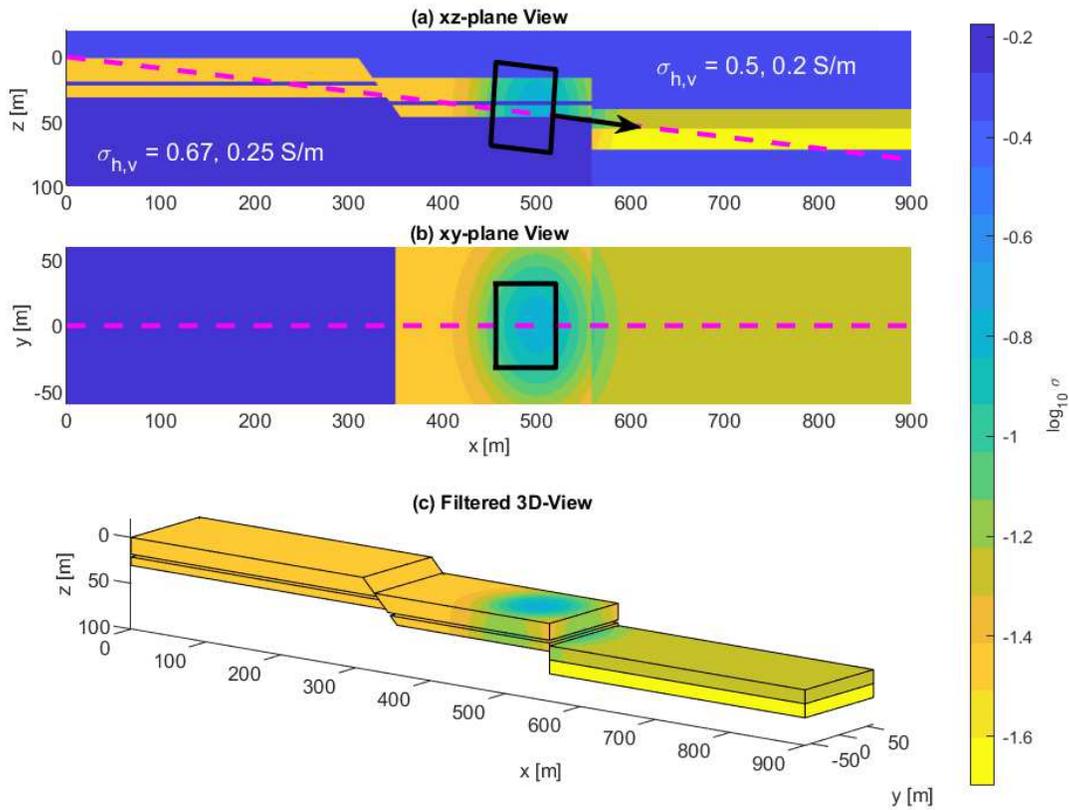}
    \caption{(a) xz-plane view at y = 0 m. (b) xy-plane view at z = 40~m. (c) 3D view of the model with the shale layers removed. The magenta dashed lines indicate the drilling trajectory. The black box is a forward modelling window example with a size of 64 $\times$ 64 $\times$ 64 m\textsuperscript{3} at one logging position.}
    \label{fig:8_ex3_model}
\end{figure}

\begin{figure}
    \centering
    \includegraphics[width=7cm,height=2cm]{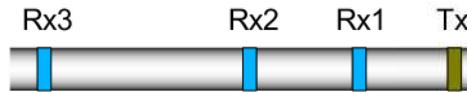}
    \caption{Illustration of an induction tool with a single transmitter and three receivers. Tx and Rx stand for transmitter and receiver, respectively.}
    \label{fig:9_tools_schematic}
\end{figure}

\begin{figure}
    \centering
    \includegraphics[width=12cm,height=6cm]{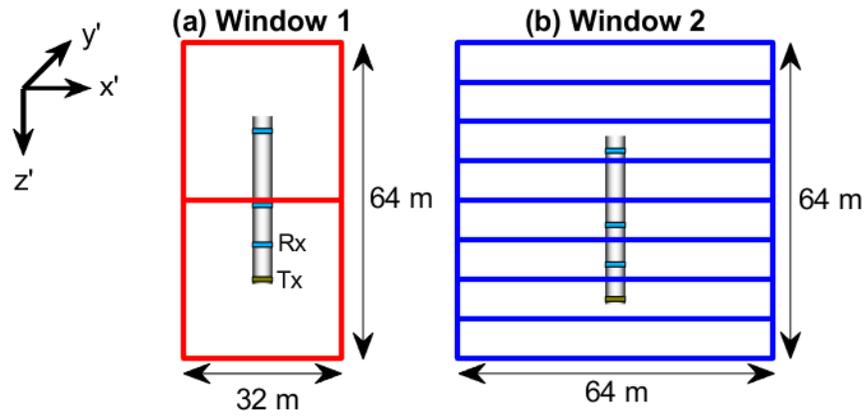}
    \caption{Domain decomposition illustration of (a) window 1 and (b) window 2. z'-axis is the drilling direction. The range in the y' direction is equal to the range in the x' direction.}
    \label{fig:10_windows_sketch}
\end{figure}

\begin{figure}
    \centering
    \includegraphics[width=7cm,height=8.5cm]{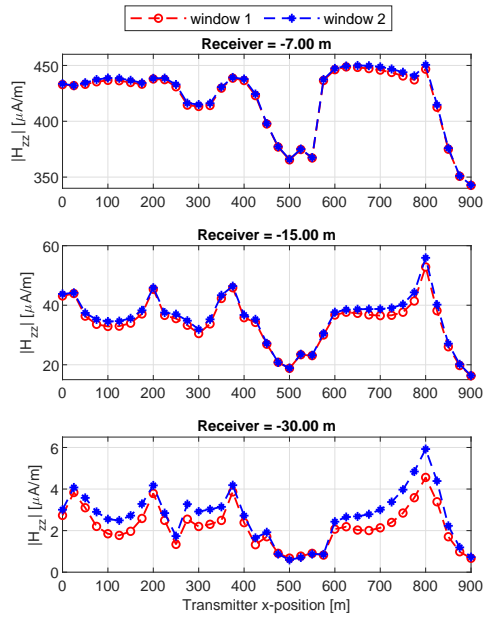}
    \caption{Measured magnitude of the z-component \textbar{}\(H_{zz}\)\textbar{} at different receiver positions across the formation.}
    \label{fig:11_ex3_log_response}
\end{figure}

\renewcommand{\thefigure}{\ref{appendix_verification}\arabic{figure}}
\setcounter{figure}{0}

\begin{figure}
    \centering
    \includegraphics{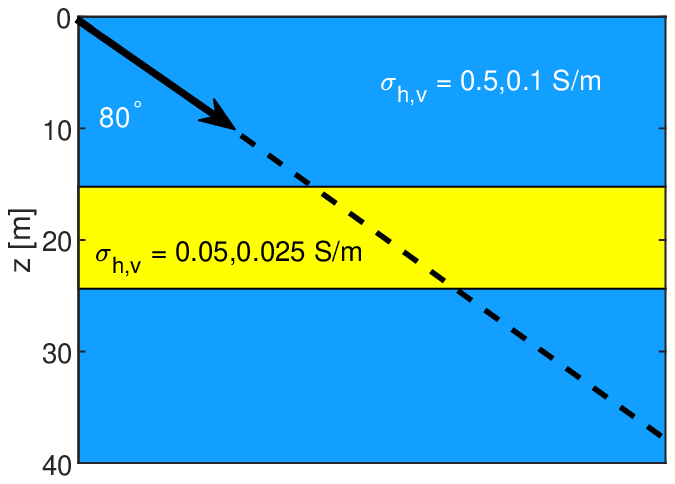}
    \caption{Sketch of logging while drilling across an anisotropic layered medium.}
    \label{fig:appendix1_model}
\end{figure}

\begin{figure}
    \centering
    \includegraphics{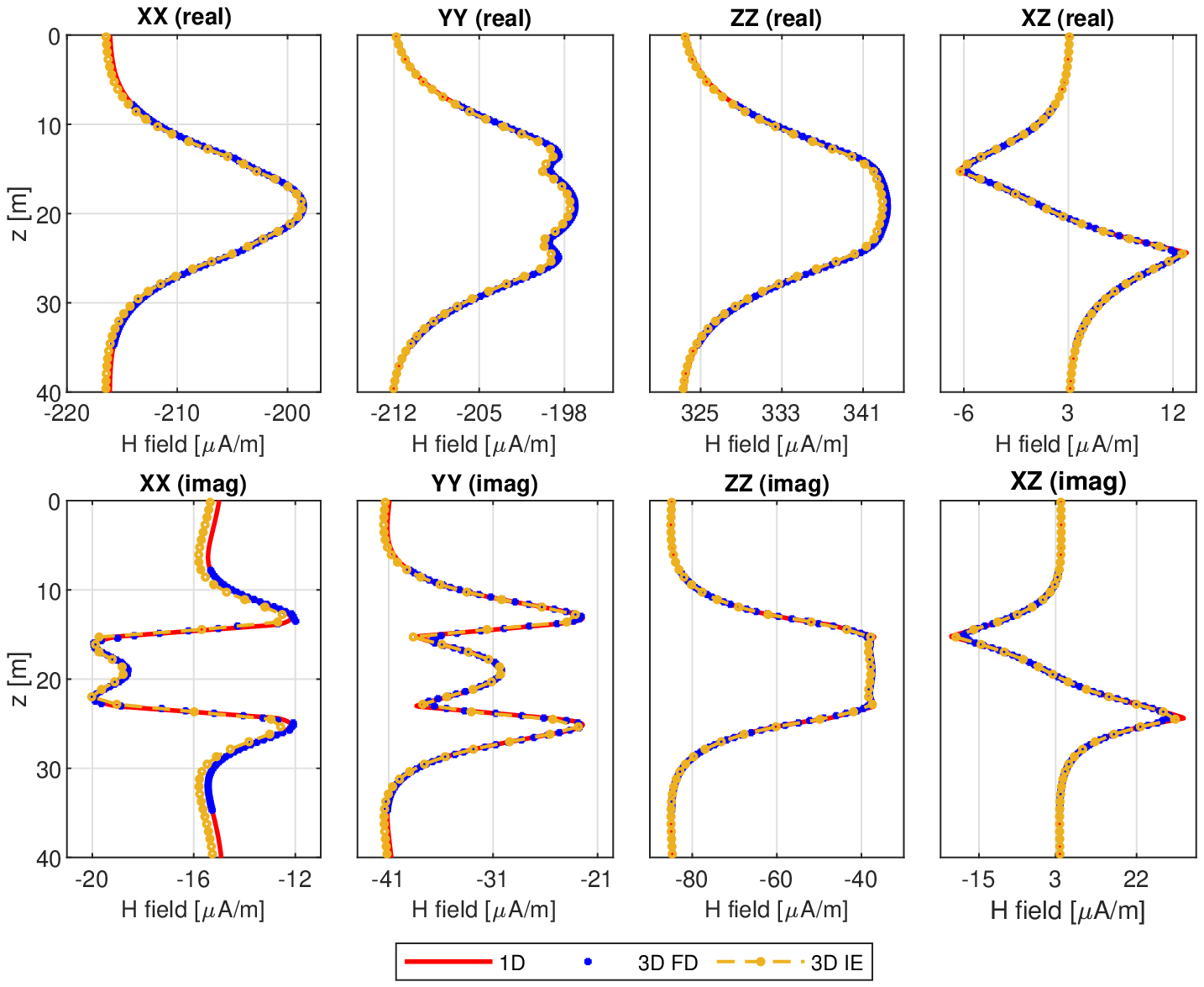}
    \caption{Comparison of the calculated magnetic field couplings with different numerical methods.}
    \label{fig:appendix2_comparison}
\end{figure}

\renewcommand{\thetable}{\arabic{table}}
\FloatBarrier
\begin{table*}
\begin{minipage}{160mm}
\caption{Performance of IE with different schemes.}
\label{table:schemes_performance}
\begin{tabular}{lllll}
  &  & Total GMRES & Total Outer &  \\
Schemes & Target $e_{sub}$ & Iterations & Iterations & Time (s) \\
\hline
\begin{minipage}{30mm}
    Full Domain IE
\end{minipage}
 & $10^{-6}$ & 97 & - & 21.82 \\
 IE-DD-GS-F
 & $10^{-6}$ & 206 & 5 & 18.04 \\
 IE-DD-GS-A
 & $e_{full}$/10 & 124 & 6 & 14.18 \\
  IE-DD-Jacobi-A
 & $e_{full}$/10 & 183 & 8 & 19.49 \\
\hline
\end{tabular}

\medskip
$e_{sub}$ and $e_{full}$ are the relative residual of the sub- and full-domain, respectively.
\end{minipage}
\end{table*}

\begin{table*}
\begin{minipage}{160mm}
\label{alg:algorithm_1}
\begin{tabular}{l}
\hline
\textbf{Algorithm 1.} IE method with Domain Decomposition Preconditioning.\\
\hline
\(\textbf{set}\ : \epsilon = \text{threshold value},~M = \text{number of sub-domains},~\) \\
\(\textbf{initialize} : \mitbf{E}^{0}:=\mitbf{E}^{(0)}, k:=1 \)\\

\(e = \frac{\left\| \mitbf{E}^{(0)} - \left( \mitbf{I} - \mitbf{\mathcal{G}}\Delta\mitbf{\sigma} \right)\mitbf{E}^{k-1} \right\|}{\left\| \mitbf{E}^{(0)} \right\|}\)\\

\textbf{while} $e >$ $\epsilon$\\
\hspace{0.25cm}\textbf{for} $i = 1:M$ \\
\hspace{0.5cm}\textbf{if} Gauss-Seidel preconditioning\\
\hspace{0.75cm}\(\mitbf{b} = \mitbf{E}^{(i,0)} + \sum_{j = 1}^{i - 1}{\mathcal{G}^{(ij)}\Delta\mitbf{\sigma}^{(j)}\mitbf{E}^{(j),k + 1}} + \sum_{j = i + 1}^{M}{\mathcal{G}^{(ij)}\Delta\mitbf{\sigma}^{(j)}\mitbf{E}^{(j),k}}\)\\
\hspace{0.5cm}\textbf{else if} Jacobi preconditioning\\
\hspace{0.75cm}\(\mitbf{b} = \mitbf{E}^{(i,0)} + \sum_{j = 1}^{M}{\mitbf{\mathcal{G}}^{(ij)}\Delta\mitbf{\sigma}^{(j)}\mitbf{E}^{(j),k}}\)\\
\hspace{0.5cm}\textbf{end if}\\
\hspace{0.5cm}\textbf{set} : initial guess = \(\mitbf{E}^{(i),k}\),~threshold = $e$/10\\
\hspace{0.5cm}$\mitbf{E}^{(i),k + 1}$ = GMRES$\left[\mitbf{A} = \left( \mitbf{I} - \mitbf{\mathcal{G}}^{(ii)}\Delta\mitbf{\sigma}^{(i)} \right), \mitbf{b}, \text{initial guess}, \text{threshold}\right]$\\
\hspace{0.25cm}\textbf{end for}\\
\hspace{0.25cm}$k = k+1$\\
\textbf{end while}\\
\hline
\end{tabular}
\end{minipage}
\end{table*}

\newpage
\FloatBarrier

\end{document}